\documentclass{article}

\usepackage{arxiv}

\usepackage[utf8]{inputenc} 
\usepackage[T1]{fontenc}    
\usepackage{hyperref}       
\usepackage{url}            
\usepackage{booktabs}       
\usepackage{amsfonts}       
\usepackage{nicefrac}       
\usepackage{microtype}      
\usepackage{lipsum}
\usepackage{graphicx}
\graphicspath{ {./images/} }

\usepackage{algorithm}
\usepackage{algorithmic}
\usepackage{calrsfs}
\DeclareMathAlphabet{\pazocal}{OMS}{zplm}{m}{n}

\usepackage{amsmath}

\title{AccEq-DRT: Planning Demand-Responsive Transit to reduce inequality of accessibility}

\author{%
  \textbf{Duo Wang}\\
  Samovar, Télécom SudParis, Institut Polytechinque de Paris\\
  19 place Marguerite Perey, 91120 Palaiseau, France \\
  duo.wang@ip-paris.fr \\
  \And
  \textbf{Andrea Araldo}\\
  Samovar, Télécom SudParis, Institut Polytechinque de Paris\\
  19 place Marguerite Perey, 91120 Palaiseau, France \\
  andrea.araldo@telecom-sudparis.eu \\
  \And
  \textbf{Mounîm A. El Yacoubi}\\
  Samovar, Télécom SudParis, Institut Polytechinque de Paris\\
  19 place Marguerite Perey, 91120 Palaiseau, France \\
  mounim.el$\_$yacoubi@telecom-sudparis.eu 
}

\begin{document}
\maketitle
\begin{abstract}
Accessibility measures how well a location is connected to surrounding opportunities. We focus on accessibility provided by Public Transit (PT). There is an evident inequality in the distribution of accessibility between city centers or close to main transportation corridors and suburbs. In the latter, poor PT service leads to a chronic car-dependency. Demand-Responsive Transit (DRT) is better suited for low-density areas than conventional fixed-route PT. However, its potential to tackle accessibility inequality has not yet been exploited. On the contrary, planning DRT without care to inequality (as in the methods proposed so far) can further improve the accessibility gap in urban areas.

To the best of our knowledge this paper is the first to propose a DRT planning strategy, which we call AccEq-DRT, aimed at reducing accessibility inequality, while ensuring overall efficiency. To this aim, we combine a graph representation of conventional PT and a Continuous Approximation (CA) model of DRT. The two are combined in the same multi-layer graph, on which we compute accessibility. We then devise a scoring function to estimate the need of each area for an improvement, appropriately weighting population density and accessibility. Finally, we provide a bilevel optimization method, where the upper level is a heuristic to allocate DRT buses, guided by the scoring function, and the lower level performs traffic assignment. Numerical results in a simplified model of Montreal show that inequality, measured with the Atkinson index, is reduced by up to 34\%.
\end{abstract}


\section{Introduction}
The accessibility of a location measures the ease for its residents to reach surrounding opportunities (e.g., schools, shops, or jobs). Here we study PT accessibility, i.e., how many opportunities per unit of time one can reach by using Public Transit (PT). Accessibility is unequally distributed in urban regions~\cite{Badeanlou2022}: it is generally good in the city centers and close to main PT corridor, and poor in the suburbs, which are thus car-dependent. This makes cities environmentally, economically and socially unsustainable~\cite{giuffrida2017investigating,Boussauw2022}.
Poor accessibility in the suburbs is a chronic problem of conventional PT. Offering a good service would require an extremely large number of lines, stops and frequencies. Considering the low population density, the cost pro-capita would be prohibitive. 

By dynamically building routes that adapt to the observed user requests, Demand-Responsive Transit (DRT) is better suited for low density areas~\cite{quadrifoglio2009methodology}. DRT is thus a promising complement to conventional PT~\cite{calabro2023adaptive}. Vehicle automation is expected to reduce the cost, thus making DRT a feasible option for PT operators~\cite{fielbaum2020strategic}. 
The literature has proposed planning strategies have aimed to minimize generalized cost, including both user and agency related cost, and found that DRT tilts the cost balance in favor of users~\cite{MarcoNie2019}. However, we believe that studying DRT only under the lenses of generalized cost minimization risks downsizing its potential benefits. We believe one of the most important benefit of DRT is its capability to improve accessibility in locations in which it would not be possible via conventional PT.

The main contribution of this paper is AccEq-DRT, an optimization method to support DRT planning, at strategic level, in an urban conurbation. The method suggests in which areas to deploy DRT and with how many vehicles. The main novelty is that the goal of our optimization is to reduce accessibility inequality, instead of traditional overall cost-minimization. Three aspects make our goal particularly challenging. \emph{First}, accessibility is the result of intricate and latent dependencies between different parts of the PT network. Therefore, a local change may affect accessibility of remote areas. Such dependencies are coupled in a non-trivial with the topology of the PT network. \emph{Second}, there is a complex demand-supply dependence, since the performance of DRT impacts the number of users that choose to use it, which in turn impacts the performance. \emph{Third}, when considering accessibility inequality one should not  only examine the geographical distribution of accessibility, but also population density. This is important in order to prevent DRT planning decisions from being heavily driven by few people who have chosen to live in isolated areas. Population and accessibility distributions may be contrasting. 

A second contribution is methodological and concerns the modeling approach. We combine a Continuous Approximation (CA) model of DRT and a graph model of conventional PT. CA allows to efficiently estimate the average performance indicators of DRT. On the other hand, the graph model of conventional PT  represents the topology of the network that already exists in a city, thus capturing the geographical dependencies in the accessibility computation. This makes the findings obtained with our model relevant for the city under study. We obtain a single multilayer graph describing a multimodal PT, containing both DRT and conventional PT, in which accessibility can be computed.

To solve the demand-supply circular dependency, our optimization procedure is bilevel: the upper level decides how many DRT vehicles to deploy in each area, while the lower level iteratively solves a transit assignment problem, to estimate how travelers distribute in the multimodal network.

Numerical results in a scenario representing Montreal show that AccEq-DRT effectively reduces accessibility inequality, with a limited fleet of buses, while also increasing overall accessibility. We release our code to run our strategy and reproduce the analysis as open source~\cite{DuoCode}.

\subsection{RELATED WORK}
Within the vast literature on Demand-Responsive Transit (DRT), the most relevant work for our paper is that concerning the integration of DRT into conventional PT. Such integrated system has been called \emph{Adaptive Transit}~\cite{calabro2023adaptive}, to emphasize its ability to exploit both the flexibility of DRT and the efficiency of conventional PT to better adapt to the spatially and temporally heterogeneous demand of big urban regions. Detailed agent-based simulation~\cite{articleBasu,Chouaki2021,Seshadri2023} or graph models~\cite{Pavone2019} has been used to evaluate system- and user-level performance of DRT or ride-sharing.

Strategic decisions for planning DRT together with conventional PT have been attracted recent interest in the community~\cite{calabro2023adaptive,MarcoNie2019,steiner2020strategic,Mahmassani2020}. The objective of the optimization in the aforementioned work is always the minimization of some generalized cost. We adopt the radically different objective of planning DRT for reducing accessibility inequality.

Improvement in accessibility obtained via DRT has been recently estimated via simulation~\cite{Nahmias2021,Bischoff2023} and geostatistical analysis~\cite{Diepolder2023}. However, we are the first to set the reduction of accessibility distribution inequality as main optimization goal of the optimization of DRT planning. 
Recent paper considering accessibility and inequality in optimization problems within transportation planning are~\cite{article5}, which added a constraint concerning equality in a classic cost minimization problem, \cite{aaaaaaaarticle6}, which proposed timetables and ticketing to favor low-income passengers and~\cite{Dai2022}, optimizing headway and stop spacing of a single line considering equality. However, none of these papers consider DRT. 

As for the methodological novelty, we note that optimization of the design of Adaptive Transit is often tackled with fully analytical models, e.g., continuous approximation~\cite{calabro2023adaptive,MarcoNie2019} or similar~\cite{fielbaum2020strategic}. Their main limit is that the PT network is represented with a regular geometrical pattern. Such idealized structures can loosely describe ``any city'', but at the end they describe ``no city'', as the pre-existent PT network in a city is always much more complex than regular geometrical shapes. In our work we focus on accessibility, which results from an intimate interplay between the actual structure of PT deployed in a city and how population distributes around it. Therefore, idealizing PT with simple shapes would lead to findings with no real relevance. This rules out the aforementioned fully analytical approaches. On the opposite side of the methodological spectrum, building a detailed simulation of the entire multimodal PT as in~\cite{Mahmassani2020} is very complex to develop and heavy to run together with multiple optimization iterations. This would overshoot the needs of preliminary findings at a strategic level. Moreover, since we are targeting new deployments of DRT services that are not currently operating, we must accept that even detailed simulation will not match the performance of such services in reality. For all these reasons, we decide to use CA to model DRT. We did not find in the literature other models combining CA for the first and last mile DRT and a graph for representing with some realism the structure of conventional PT in a certain city. A similar work in this sense is~\cite{Wong2003}, which is however dealing with completely different dynamics (car flow dynamics, where a graph represents highways and CA models represent local streets). Moreover, we optimize DRT deployment, while~\cite{Wong2003} does not involve any strategic decision. Indeed, most models

\section{MODEL}

\begin{figure}
    \centering
    \includegraphics[width=0.6\textwidth]{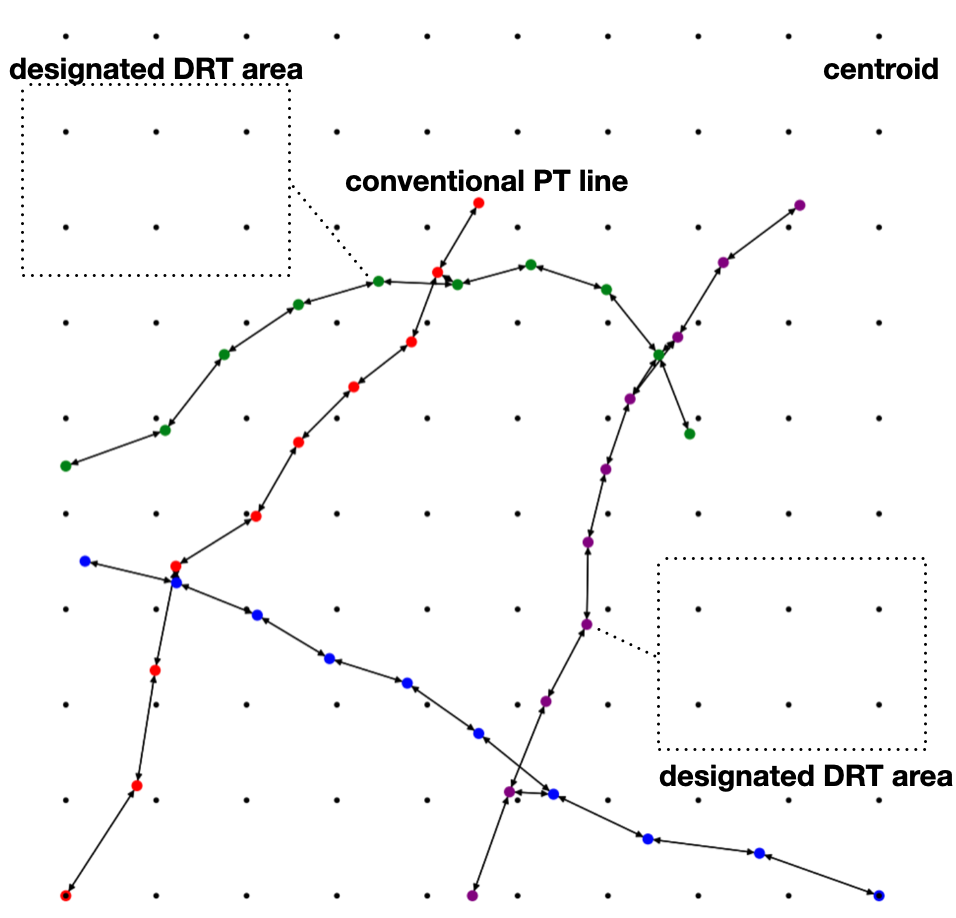}
    \caption{Multimodal Public Transit} 
    \label{fig:passenger-breakdown}
\end{figure}

We study a multimodal urban transit system, as in Fig.~\ref{fig:passenger-breakdown}, composed of~\cite{calabro2023adaptive}:
\begin{enumerate}
    \item Conventional PT, consisting of metro, train or bus lines, with regular routes and schedules.
    \item Demand-Responsive Transit (DRT), provided by bus, acting as a feeder for conventional PT and serving the First Mile and Last Mile.
\end{enumerate}

\begin{figure}
    \centering
    \includegraphics[width=1.0\textwidth]{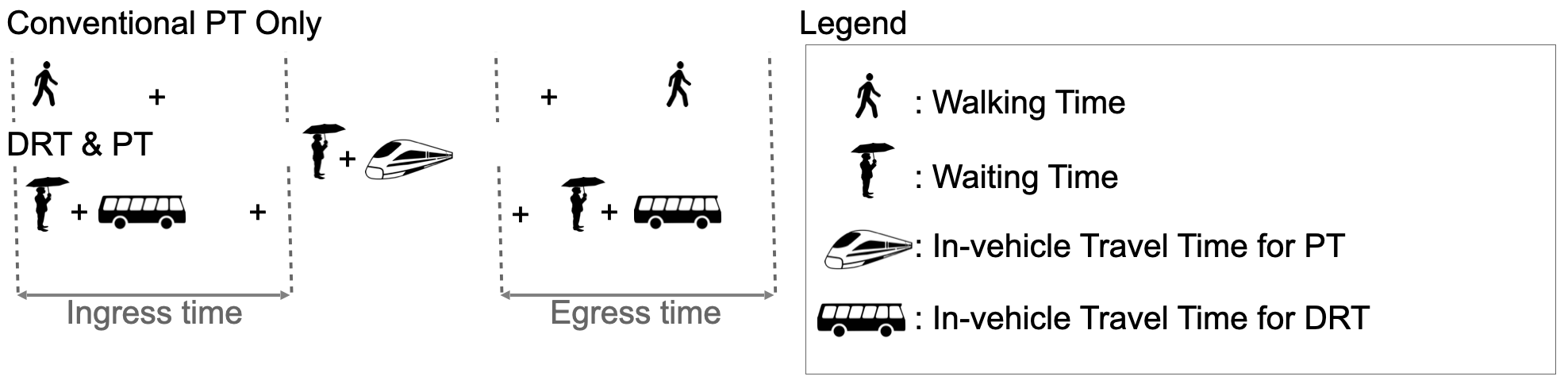}
    \caption{Components of the access, egress and waiting time for the case with conventional PT only and the case with DRT.} 
    \label{fig:passenger-1}
\end{figure}
Fig.~\ref{fig:passenger-1} represents the types of trips a user can perform in the considered transit system. 

We now describe the model that we use as a base for AccEq-DRT optimization. After explaining the tessellation of the study area, we first explain how we model conventional PT and DRT separately and how we merge such models into a single multilayer graph. We then explain how we generate the demand.

\subsection{Tessellation of the study area}
We partition the study area with a regular tessellation. In the numerical results we adopt square tiles of 1~Km\textsuperscript{2}, but any regular shape can be employed. The center of each tile is called \emph{centroid}. We denote by $\pazocal{C} := \{c_1,\dots,c_{M}\}$ the set of all \emph{centroids}. By slight abuse of notation, we will sometimes refer to ``tile $c_i$'' to denote the ``tile around centroid $c_i$''. We will thus use ``tile'' and ``centroid'' interchangeably, when this does not create ambiguity. \emph{Opportunities} are distributed in the tiles. They can be shops, jobs, schools. In the numerical results, for simplicity, we will consider people as opportunities. We assume tiles are small enough so that the distance between any point in a tile and the corresponding centroid in the center is negligible. We can thus approximately assume all trips generate from a centroid and end at another centroid. Note that, as in several studies on accessibility~\cite{Loreto2019}, we prefer regular tassellation to standard zoning. This will allow simple development of the CA model (Equations~\ref{eq:n_passengers}-\ref{eq:ingress}).

\subsection{Graph model of conventional PT}
\label{sec:Model1}

We model conventional PT as a graph $\pazocal{G}$ (see Figure~\ref{fig:passenger-breakdown}). $\pazocal{G}$ is composed of multiple lines. A line $l$ has a headway $t_l$, which is the time between two vehicles in the same direction. We approximate the average waiting time for line $l$ as~$w_l=t_l/2$ as in~\cite{Mahmassani2020}. Line $l$ is a sequence of stations, $s_1,\dots,s_{K_l}$, linked by the edges. Let us denote with $t(s_j, s_{j+1})$ the time spent by a vehicle to go from $s_j$ to $s_{j+1}$.To allow boarding and alighting at station $s_{j}$, the vehicle stops for a dwell time $t_{s_j}$. 
If a passenger is at a transfer station $s_{j}$ and changes from PT line $l$ to line $l'$, we consider an average waiting time $t^{l'}$ at station $s_{j}$.

We also include in $\pazocal{G}$ the set of centroids and edges between any centroid and all the others and betweren any centroid and all stations. 

For a trip from centroid $c_i$ and $c_j$ , the traveler can choose between different modes of travel. For example, a traveler could perform the entire trip by walking at speed $v_{walking}$ or they could walk from the origin centroid $c_i$ to a conventional PT station, go via conventional PT to another station, and from there walk to the destination centroid $c_j$. A traveler always selects the shortest time path between a centroid and another (for brevity, we will call it ``shortest path''). This is more accurate than~\cite{MarcoNie2019,calabro2023adaptive} in which shortest distance path is selected. 

Note that the graph $\pazocal{G}$ represents a PT network in a certain time slot, e.g. 1 hour, within which we can assume that line routes and headway values do not change. To obtain a planning over a day, the optimization we will present here should be repeated in the different time-slots.

\subsection{Continuous Approximation (CA) model of DRT}

\begin{figure}
    \centering
    \includegraphics[width=1.0\textwidth]{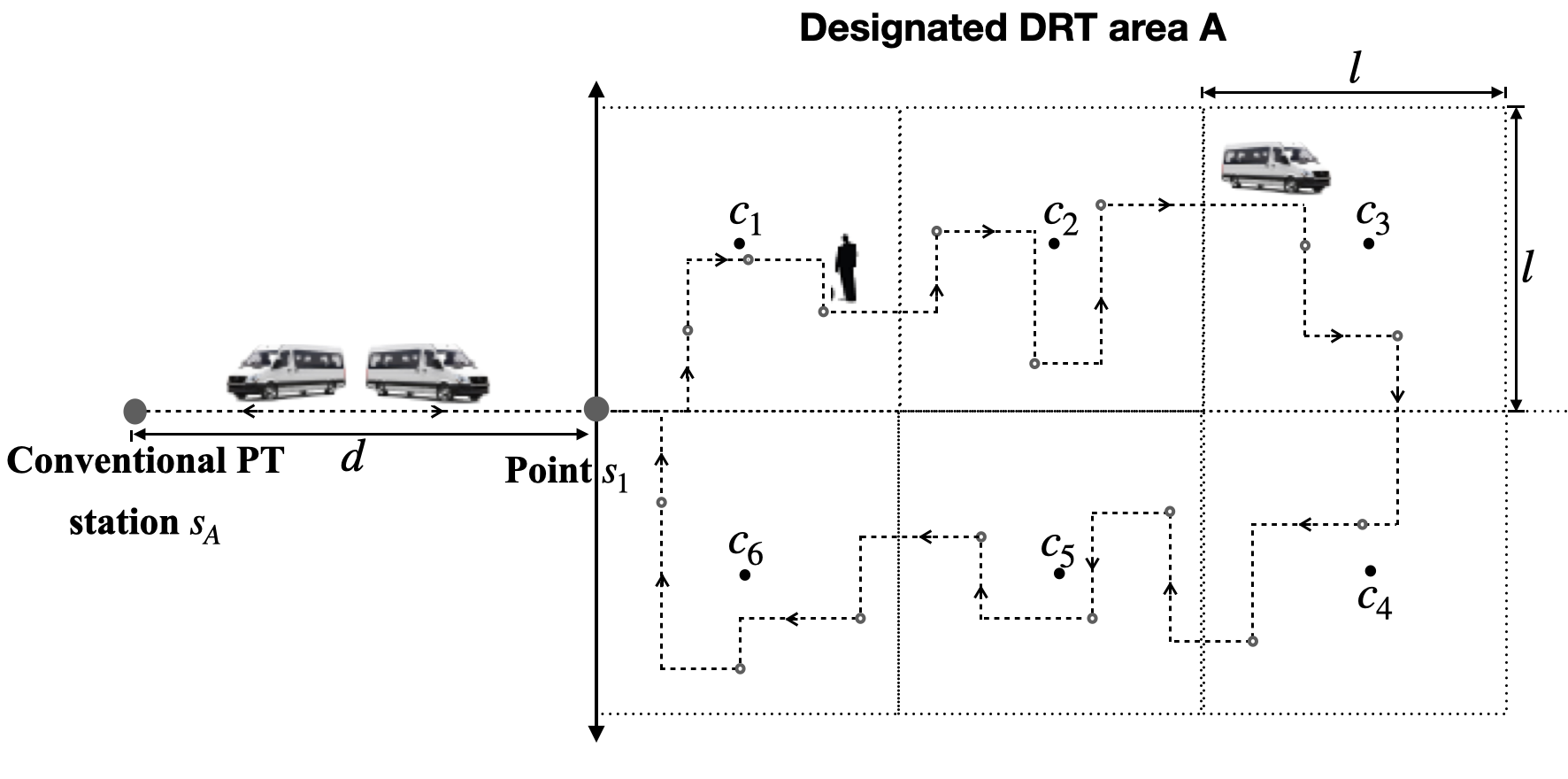}
    \caption{Model of Demand Response Transit (DRT).} 
    \label{fig:Fig3}
\end{figure}

Due to financial constraints, there is a limitation on the total number of buses that can be deployed for DRT services. Because of this, only selected \emph{designated DRT areas} can be served via DRT. Vehicles will depart from a selected conventional PT station at regular times to pick up travelers in the designated DRT area. As a result, passengers in the service area will travel as shown in Figure~\ref{fig:passenger-1}. A passenger may wait until a DRT bus arrives (waiting time) and spend some time in the bus (in-vehicle travel-time) to the conventional PT station. Symmetrically, to reach the final destination from an egress PT station, a passenger must walk or use another DRT service, if available. 

The DRT model is shown in Figure~\ref{fig:Fig3}. For a generic designated DRT area $A$ associated to conventional PT station $s_A$, a DRT bus departs from $s_A$ and proceeds to point $s_1$, without performing any pick up/drop off between $s_0$ and $s_1$. Starting from $s_1$, the bus serves each tile in sequence and then return to station $s_A$ via point $s_1$. 

Our DRT model is based on~\cite{quadrifoglio2009methodology} and Sec. 4.4 of~\cite{calabro2023adaptive} and relies on similar assumptions. A DRT bus can only travel along the horizontal and the vertical direction following a no-backtracking policy~\cite{quadrifoglio2009methodology}. The DRT service has a headway $h(A)$, and we assume bus capacity is enough to accommodate all passengers during headway $h(A)$. Demand is uniformly distributed within each tile. Let $n(A)$ denote the number of pickups and dropoffs requested within area $A$ between the departure of a DRT bus and the next. 
\label{ln:originating} To compute $n(A)$, let~$\phi_i^\text{DRT, out}$ denote  the number of passengers per unit of time that originate in tile $c_i$ and that use DRT for their first mile. Let $\phi_i^\text{DRT,in}$ denote the number of passengers per unit of time having tile $c_i$ as destination and using DRT for their last mile. If $l$ is the length of a tile, the number of pickups and dropoffs in designated DRT area $A$ is the following sum over all centroids $c_i$ contained in designated DRT area~$A$:

\begin{align}
&& n(A) &= h(A)\cdot l^2\sum_{c_i \in A}  \phi_i^\text{DRT, out} 
            + 
        h(A)\cdot l^2\sum_{c_i \in A} \phi_i^\text{DRT, in}.
\label{eq:n_passengers}
\end{align}

The two terms in the above sum correspond to the number of pickups and dropboffs, respectively.
The cycle length measures how many Km the DRT bus travels from the moment it leaves $s_A$ to the moment it comes back. The expected value of cycle length $CL_{DRT}(A)$ is computed precisely as Eqn. 10 of~\cite{quadrifoglio2009methodology}, with length $L = \frac{Kl}{2}$ and width $W = 2l$:

\begin{align}
&& CL_{DRT}(A) &= 2d + Kl \frac{n(A)}{n(A)+1} + n(A)\cdot \frac{l}{3}  +\frac{4l}{3}
\label{eq:cycle_length}
\end{align}

The expected time required to complete a cycle is: 

\begin{align}
&& C_{DRT}(A) &= \frac{CL_{DRT}(A)}{v_{DRT}} + \tau_s  \cdot n(A) + \tau_T,
\label{eq:cycle_time}
\end{align}

\noindent where $v_{DRT}$ the DRT bus's speed, $\tau_s$ the time lost per stop (assuming a single passenger served per stop) and $\tau_T$ the terminal dwell time~\cite{calabro2023adaptive} of the DRT bus and the terminal $s_A$. If we assign $N_A$ buses to this designated DRT area,  headway $h(A)$ is

\begin{align}
&& h(A) &= C_{DRT}(A) / N_A
\label{eq:hcl}
\end{align}

\noindent We now compute ingress time~$T_{in}(c_k)$ at centroid $c_k\in A$, i.e., the time it would take for a user at tile $c_k$ from the moment they requested a trip to the moment in which they arrive at the designated station $s_A$. Symmetrically, we define the egress time. Assume designated DRT area $A$ has $2 \times \frac{K}{2} $ tiles (K=6 in Figure~\ref{fig:Fig3}). For centroid $c_k\in A$ let $n(c_k)$ denote the number of pickups and dropoffs requested within tile $c_k$ between the departure of a DRT bus and the next. The expected value of ingress and egress access time for centroid $c_k$ is: 

\begin{align}
&& T_{in}(c_k) &= T_{out}(c_k) = \frac{h(A)}{2} + \frac{ \frac{n(c_k)}{2} + \sum_{i=k+1}^{K}n(c_i)}{n(A)} \cdot \left(C_{DRT}(A) - \frac{2d}{v_{DRT}} \right) + \frac{d}{v_{DRT}}
\label{eq:ingress}
\end{align}

\noindent The three components of Equation~\ref{eq:ingress} calculate the average waiting time for the DRT bus, the average in-vehicle time from $c_k$ to $s_1$, and the average in-vehicle time from $s_1$ to $s_A$, respectively (as shown in Figures~\ref{fig:passenger-1} and~\ref{fig:Fig3}). We now explain the calculation method of the second component. Because the pickup and dropoff density $n(c_i)$ within different tiles $c_i$ in the considered designated DRT area is not the same, the calculation formula differs from~\cite{quadrifoglio2009methodology}. The fraction of pickups and dropoffs demand that are in tile $c_k$ is $\frac{n(c_k)}{n(A)}$. Therefore $\frac{n(c_k)}{n(A)} \cdot  \left(C_{DRT}(s_0,s_1) - \frac{2d}{v_{DRT}}\right)$ is the average time it takes to pass through this tile. A passenger who gets on the bus at $c_k$ needs on average to go through half of tile $c_k$ and all the subsequent tiles, which amount to fraction $\frac{\sum_{i=k+1}^{K}n(c_i)}{n(A)}$ of pickups and dropoffs (Fig.~\ref{fig:Fig3}).

\subsection{Multilayer graph}

We denote the multilayer graph of discrete conventional PT and continuous approximation of DRT by $\pazocal{M}$ as shown in Figure~\ref{fig:passenger-breakdown}. 
In addition to all nodes and edges of $\pazocal{G}$, $\pazocal{M}$ also includes \emph{DRT edges}. These edges are between any centroid in any designated DRT area $A$ and the corresponding conventional PT station $s_A$ (Fig.~\ref{fig:Fig3}). Their values represent the expected ingress and egress times (Eqn.~\ref{eq:ingress}).
By optimizing the number $N_A$ of buses deployed in designated DRT area $A$, we change the times on such DRT edges.

\subsection{Demand model}
To use Eqn.~\ref{eq:n_passengers} and below, we need to know $\phi_i^\text{DRT, out}$ and $\phi_i^\text{DRT, in}$. We now explain how we obtain them for a certain time-slot. Let $\rho_i$ be the population density in tile~$c_i$ and $\textit{Trip}$ the trip rate, i.e., the percentage of residents performing any trip in that time-slot.
The trip density $\phi_i$ originating in centroid $c_i$ is
\begin{align}
    && \phi_i &= \textit{Trip} \cdot \rho_i.
    \label{eq:phi_i}
\end{align}

Let $\sigma_j$ be the amount of opportunities in tile $c_j$. Note that $\sigma_j$ represents the ``attractiveness'' of $c_j$. Adopting a singly constrained gravity model (formulas (1) and (2) in~\cite{ABDELAAL2014677}), the trip density $\phi_{ij}$ from centroid $c_i$ to $c_j$ is :

\begin{align}
&& \phi_{ij} &= \phi_i \cdot \frac{ \sigma_j \cdot e^{-\beta T(c_i,c_j)}}{ \sum_{k=1}^{M} (\sigma_k \cdot e^{-\beta T(c_i,c_k) } )} \,\,\,\,\,[trips/hour]
\label{eq:trips-2}
\end{align} 
being $\beta$ a constant called dispersion parameter and $T(c_i,c_j)$ the shortest time to go from $c_i$ to $c_j$ on the multimodal graph $\mathcal{M}$. Equation~\ref{eq:trips-2} captures the fact that travelers tend to go to places with many opportunities and reachable in a short time. Observe that, in order to go from $c_i$ to $c_j$, the shortest path might use DRT in the first mile, or in the last mile, or in both, or in none of them. To account for this, let us assume that $c_i$ is inside a designated DRT area $A$ and that $s_A$ is the associated conventional PT station (Fig.~\ref{fig:passenger-breakdown}). We introduce the indicator function $\mathbb{I}_{(c_k, c_j)}^{s_A,\text{first}}$, which is $1$ if the shortest path from $c_k$ to $c_j$ passes by $s_A$ (a passenger would enter conventional PT via $s_A$). Symmetrically, indicator function~$\mathbb{I}_{(c_j, c_k)}^{s_A,\text{last}}$ is $1$ if the shortest path from $c_j$ to $c_k$ passes by $s_A$ (a passenger would exit conventional PT from $s_A$ to reach $c_k$ via DRT). Otherwise, such indicators are $0$. The numbers of or DRT users originating and directed to centroid $c_k$ are, respectively (see definitions in page~\pageref{ln:originating})

\begin{align}
&&\phi_k^\text{DRT, out}&=\sum_{c_j\notin A} \phi_{k,j} \cdot \mathbb{I}_{(c_k, c_j)}^{s_A,\text{first}}
\label{eq:flow-out}
\\
&&\phi_k^\text{DRT, in}&=\sum_{c_j\notin A} \phi_{j,k} \cdot \mathbb{I}_{(c_k, c_j)}^{s_A, \text{last}}
\label{eq:flow-in}
\end{align}

\section{DRT planning strategy}

We first introduce accessibility, which is at the base of our approach, and the inequality indicator we will use to evaluate the quality of our method. We then describe our optimization procedure. 

\subsection{Accessibility and inequality}
\label{sec:Accessibility}

The accessibility of a centroid measures how well it is connected to the surrounding opportunities (schools, people, businesses). We select the classic \emph{gravity-based} definition of accessibility~\cite{/content/paper/8687d1db-en}, which can be intended as the number of opportunities reachable in the unit of time. The accessibility of centroid $c_i\in \pazocal{C}$ is
\begin{align}
&& acc(c_i) &= \sum_{c_j} \frac{\sigma_j}{T(c_i,c_j)}
\label{eq:acc}
\end{align}
where $ T(c_i,c_j) $ is the shortest time to go from $c_i$ to $c_j$ and $ X(c_j) $ is the amount of opportunities of the tile having centroid $c_j$. In broad terms, a large $acc(c_i)$ indicates that the centroid $c_i$ is well connected to opportunities. When optimizing for equality, we will focus in particular on the ``most unfortunate'' individuals, i.e., the ones with the lowest accessibility. For all individuals $v$ resident in tile $c_i$, we define their accessibility as
\begin{align}
&& acc(v) &= acc(c_i)
\label{eq:acc_1}
\end{align}

We denote by~$\pazocal{V}^{m\%}$ the set of the unfortunate $m$\% people. We define the $m$ accessibility of a graph $\pazocal{G}$ as:
\begin{align}
&& Acc(\pazocal{G};m) &=  \frac{1}{|\pazocal{V}^{m\%}|}\sum_{v \in \pazocal{V}^{m\%}}acc(v),
\label{eq:Acc-graph}
\end{align}
where $|\cdot|$ denotes the cardinality of a set.
With $m=100$, $\pazocal{V}^{100\%}=\pazocal{V}$, where $\pazocal{V}$ denotes the set of all individuals. In this case, Eqn.~\ref{eq:Acc-graph} corresponds to the average individual's accessibility 
\begin{align}
    && \overline{acc}(\pazocal{G})&=Acc(\pazocal{G}, 100).    
    \label{eq:Acc-graph-avg}
\end{align}

\vspace{0.5cm}
To quantify the benefits of our AccEq-DRT strategy, we measure the inequality of the distribution of accessibility, before and after our optimization. For simplicity, we adopt in this paper indicators of horizontal inequity (also called inequality). Computing vertical inequity indicators would require to consider socio-economic characteristics of individuals~\cite{Caggiani} and would be outside of the scope of this paper. We quantify inequality via the well-known Atkinson index on the set of individuals accessibility values $\{acc(v)| \text{individual }v\in\pazocal{V}\}$. From~\cite{Costa2019inequality}, setting $\epsilon=2$, the Atkinson inequality index is:
\begin{align}
&& Atk(G) &= 1 - \frac{1}{\overline{acc}(\pazocal{G})} \cdot \left( \frac{1}{|\pazocal{V}|} \sum_{v\in\pazocal{V}} [acc(p)]^{-1} \right)^{-1}
\label{eq:atkinson}
\end{align}
where $\overline{acc}(\pazocal G)$ is the average individual accessibility (Eqn.~\ref{eq:Acc-graph-avg}).
The Atkinson index goes from $0$ (perfect equality) to $1$ (maximum inequality).

\subsection{AccEq-DRT Algorithm}
\label{sec:Greedy}

Given the limited number of buses $N$ that can be utilized for DRT and the conventional PT network $\pazocal{G}$, we propose a AccEq-DRT algorithm (Algorithm~\ref{fig:algo_2}) to allocate DRT resources efficiently, in order to reduce accessibility inequality. 

Our basic idea is to identify the regions that need large accessibility improvements and deploy DRT there. If DRT is already deployed DRT in that region, the number of DRT buses is increased. To find the regions in need of improvement, we need a score that quantifies such a need, taking into account both population and accessibility. Intuitively, the regions needing improvements are those that suffer from low accessibility and where many people live. We thus need a score function that matches this intuition in order to guide our algorithm. To build such a score, we order tiles from the least populated to the most. Let us denote $\textit{Rank}(\rho_k)$ the position of centroid $c_k$ in the aforementioned list. Separately, we also order tiles from the least accessibility to the most, and denote $\textit{Rank}( acc(c_k) )$ the position of centroid $c_k$ for accessibility rank list. For each centroid $c_k \in \pazocal{C}$, the score of centroid $c_k$:

\begin{align}
&& Score(c_k) &= \textit{Rank}(\rho_k) + \alpha  \cdot \big( |\pazocal{C}| - \textit{Rank}( acc(c_k) ) \big),
\label{eq:index_c}
\end{align}
where $\alpha$ is a constant that weights the importance given to accessibility in the score, with respect to population. If the population of $c_k$ is large and the accessibility is low, the score will be large, denoting that $c_k$ ``needs'' accessibility improvement, as it will benefit many people. Recall that DRT is not deployed within a single tile, but in designated DRT areas with $K$ tiles each, as shown in Figure~\ref{fig:Fig3}. Therefore, we need to associate a score to each potential designated DRT area $A$:
\begin{align}
&& \textit{SCORE}(A) &= \frac{1}{K} \sum_{c_k \in A} Score(c_k)
\label{eq:index}
\end{align}
A high score for designated DRT area $A$ means that large population suffers from low accessibility and that by deploying DRT there (or by increasing DRT buses) we can potentially reduce overall inequality.

Our approach is a bilevel optimization procedure. The upper level, Algorithm~\ref{fig:algo_2}, performs planning of the DRT fleet over the study area, iteratively assigning DRT buses to the designated areas with the highest score. At each iteration, the lower level, Algorithm~\ref{fig:algo_1}, is invoked to solve transit assignment. This is crucial, because increasing the number of DRT buses in a designated DRT area $A_i$ reduces ingress and egress travel times $T_\textit{in}(c_k), T_\textit{out}(c_k)$ for passengers originated or directed to any tile $c_k$ contained in $A_i$ (see Eqn.~\ref{eq:n_passengers}-\ref{eq:ingress}). However, this encourages more passengers to choose DRT, which in turn increases the aforementioned ingress and ingress travel times. Transit assignment implemented in Algorithm~\ref{fig:algo_1} is static and deterministic (travelers always deterministically choose the shortest time path), which makes our computation simple and computationally efficient. Algorithm~\ref{fig:algo_1} alternatively updates (i)~traveler flows $\phi_k^{DRT,out}, \phi_k^{DRT,in}, \phi_{k,j}, \phi_{j,k}$ entering or exiting designated DRT area $A$ and (ii)~travel timesm until convergence is observed (flows do not change more than 5\%).

\begin{algorithm}
	\caption{ AccEq-DRT Algorithm for deploying DRT services}
	\begin{algorithmic}[1]
        
		\STATE \textbf{Input} Multilayer graph $\pazocal{M}$ with stations $\pazocal{S}$ and centroids $\pazocal{C} = \{ c_1,\dots,c_M\}$.\\
        \,\,\,\,\,\,\,\,\,\,\,\,\,\,\,\ The limited number of DRT buses $N$. \\
        \,\,\,\,\,\,\,\,\,\,\,\,\,\,\,\ The set of all centroids' demand density: $\Phi=\{ \phi_1,\dots,\phi_M\}$.\\
        \STATE \textbf{Initialization} 
        \STATE Partition the study region in a set $\pazocal{A}$ of candidate designated DRT areas
        \STATE Set the center point $s_0$ in the left side of each $A\in\pazocal{A}$ (Fig.~\ref{fig:Fig3})
        \STATE Associate to each $A\in\pazocal{A}$ conventional PT stop $s_A$ that is closest to $s_0$
        \STATE Compute all travel times $T(c_k,c_j)$ and shortest paths $P(c_k, c_j)$ for all origin-destination pairs $(c_i,c_j)$ assuming no DRT.
        \STATE Compute all demand $\phi_{k,j}$ for all pairs $(c_i,c_j)$ using Eqn.~\ref{eq:trips-2}.
        \STATE \textbf{For} step $i\leftarrow 1$ to $N$:
        \STATE \,\,\,\,\,\,\,\, For each centroid $c_k \in \pazocal{C}$, calculate $Score(c_k)$ via Equation~\ref{eq:index_c}. \\
        \STATE \,\,\,\,\,\,\,\, For each area $A \in \pazocal{A}$, calculate $\textit{SCORE}(A)$ via Equation~\ref{eq:index}.\\ 
        \STATE \,\,\,\,\,\,\,\, Find the area $A_i$ with the largest $\textit{SCORE}(A)$. \\
        \STATE \,\,\,\,\,\,\,\, Update the number of DRT buses $n_{A_i} = n_{A_i} + 1 $.\\
        \STATE \,\,\,\,\,\,\,\, Update multilayer graph $\pazocal{M}$ via Algorithm~\ref{fig:algo_1}.\\
        \STATE \textbf{End For}
        \STATE \textbf{Return} Updated multilayer graph $\pazocal{M}$.
	\end{algorithmic}  
 \label{fig:algo_2}
\end{algorithm}

\begin{algorithm}
\begin{small}
	\caption{ Transit assignment within a designated DRT area $A$}
	\begin{algorithmic}[1]        
        \STATE \textbf{Input} All data used within Algorithm~\ref{fig:algo_2}.\\
        \STATE \textbf{While} \textbf{flows are not stable} (either $\sum_{k=1}^{K}\phi_k^\text{DRT, out}$ or $\sum_{k=1}^{K}\phi_k^\text{DRT, in}$ change by more than 5\%):
        \STATE \,\,\,\,\,\,\,\, \textbf{For} each tile $c_k$ inside designated DRT area $A$:
        \STATE \,\,\,\,\,\,\,\,\,\,\,\,\, Update $\phi_k^\text{DRT, out}$ via Eqn.~\ref{eq:flow-out}.\\
        \STATE \,\,\,\,\,\,\,\,\,\,\,\,\, Update $\phi_k^\text{DRT, in}$ via Eqn.~\ref{eq:flow-in}.\\
        \STATE \,\,\,\,\,\,\,\, \textbf{End For}
        \STATE \,\,\,\,\,\,\,\, Replacing $\phi_k^\text{DRT, out}$ and $\phi_k^\text{DRT, in}$ ($\forall c_k\in A$) into Eqn.~\ref{eq:n_passengers}, express $n(A)$ as a function of $h(A)$.
        \STATE \,\,\,\,\,\,\,\, Replace the above expression of $n(A)$ into the system of equations~\ref{eq:cycle_length}-\ref{eq:hcl} to find $h(A)$.
        \STATE \,\,\,\,\,\,\,\, \textbf{If} $h(A) > 0$:\\
        \STATE \,\,\,\,\,\,\,\,\,\,\, \textbf{For} each tile $c_k$ inside designated DRT area $A$:
        \STATE \,\,\,\,\,\,\,\,\,\,\,\,\,\,\,\, Compute $T_{in}(c_k)$ and $T_{out}(c_k)$ from Eqn.~\ref{eq:ingress}.\\
        \STATE \,\,\,\,\,\,\,\,\,\,\,\,\,\,\,\, Update travel times $T(c_k, c_j)$ and $T(c_j, c_k)$, $\forall c_k\notin A$, together with the \\
        \,\,\,\,\,\,\,\,\,\,\,\,\,\,\,\, respective travel times $T(c_k, c_j), T(c_j, c_k)$.\\
        \STATE \,\,\,\,\,\,\,\,\,\,\,\,\,\,\,\, Update demand density $\phi_{k,j}$ and $\phi_{j,k}$ $\forall c_j\notin A$ via gravity model (Eqn.~\ref{eq:trips-2}).\\
        \STATE \,\,\,\,\,\,\,\,\,\,\, \textbf{End For}
        \STATE \,\,\,\,\,\,\,\, \textbf{ELSE}:\\
        \STATE \,\,\,\,\,\,\,\,\,\,\,\,\, Exit the algorithm.
        \STATE \textbf{End For}
        \STATE \textbf{Return} Updated data within Algorithm~\ref{fig:algo_2}.
        
	\end{algorithmic}  
 \label{fig:algo_1}
\end{small}
\end{algorithm}

\section{Evaluation}

\subsection{Considered scenario}
\begin{figure}
    \centering
    \includegraphics[width=1.0\textwidth]{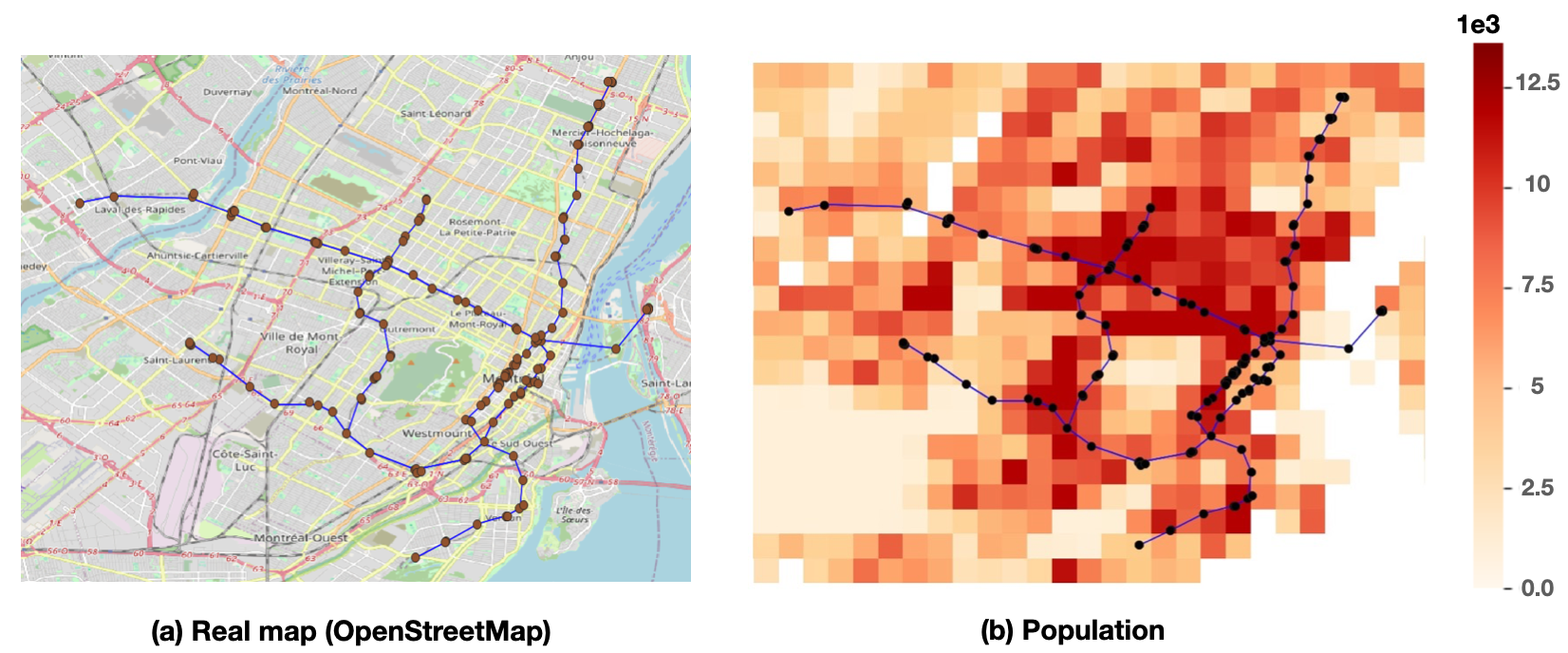}
    \caption{Metro network of Montreal and population heatmap.} 
    \label{fig:passenger-0}
\end{figure}
We build a simplified model of Montreal. Conventional PT consists of the 4 metro lines of Montreal. We do not include in our analysis fixed-line buses, for simplicity. Fixed buses are identical to metro, from the modeling point of view. So adding them would only imply more time for pre-processing the data and running the optimization, but no additional methodological complexity. Since this paper is not intended to give a detail case-study on Montreal, but to showcase the method, we decided to only consider metro lines and DRT.

From the General Transit Feed Specification (GTFS) data of Montreal~\cite{GTFS_data_montreal,population_montreal}, we take the station location, the sequence of stations of all lines, the dwell times at each station and the travel time between a station and the other. We get the headway of each metro line during peak hours from~\cite{Wiki_metro}.
In order to filter out non urbanized regions, we do not include in the study areas the centroids that are located more than 5 km far from a metro line and the tiles with $0$ population. Scenario parameters are in Table~\ref{tab:parameters}.
\begin{table}[!ht]
	\caption{Scenario parameters}\label{tab:parameters}
    \begin{small}
    \begin{center}
		\begin{tabular}{l l l l}
			Name & Parameter & Value & Reference  \\\hline
			headway &$t_l$   & 1, 2, 4 min & GTFS data~\cite{GTFS_data_montreal}  \\
			average waiting time &$w_l$ & 0.5, 1, 2 min& Equation 2 of~\cite{Mahmassani2020}  \\
            walking speed &$v_{walking}$   & 4.5 km/h & Google Maps  \\\hline
            average trip rate per peak hour &$\textit{Trip}_\text{peak}$   & 0.16 trips/person/hour & ~\cite{ABDELAAL2014677} and  ~\cite{calabro2023adaptive} \\
            total time of peak hour per day &$t_\text{peak}$   & 5.5 h & ~\cite{calabro2023adaptive}  \\
            total time of off-peak hour per day &$t_\text{off-peak}$   & 8.5 h & ~\cite{calabro2023adaptive}  \\
            dispersion parameter of Eq.~\ref{eq:trips-2} &$\beta$   & 0.12 & ~\cite{ABDELAAL2014677})  \\\hline
            length of a tile &$l$       & 1 km & ~\cite{population_montreal}  \\
            number of tiles in a DRT area &$K$       & 6 & -  \\
            time lost per DRT pickup or dropoff  &$\tau_s$       & 32 s & ~\cite{calabro2023adaptive}  \\
            terminal dwell time &$\tau_T$       & 1 min & ~\cite{calabro2023adaptive}  \\
            DRT bus speed &$v_\text{DRT}$       & 25 km/h & ~\cite{calabro2023adaptive}  \\
            \hline
            weight for the accessibility of Eqn.~\ref{eq:index_c} &$\alpha$       & 1, 2 & -  \\
            \hline
		\end{tabular}
	\end{center}
 \end{small}
\end{table}

For the sake of simplicity, we consider population as the opportunity. Therefore, our accessibility metric measures the ease to reach other people. Such accessibility metric has been called \emph{sociality score}~\cite{Loreto2019}. Therefore densities of opportunity and of population in each time are the same, i.e., $\sigma_i=\rho_i$.
\label{line:symmetric}
Under such a simplification, if a tile $c_i$ has a large population density, the trip origination density $\phi_i$ is high, but also the trip arrival density is high. We can thus adopt a symmetric trip assumption, which allows to simplify, similar to~\cite{calabro2023adaptive}, the average number of pickups and dropoffs~\eqref{eq:n_passengers} as: 
\begin{align}
&& n(A) &= 2h(A)\cdot l^2\sum_{c_i \in A}^{}\phi_i
\label{eq:n_passengers-simple}
\end{align}

\subsection{Demand generation}
Let $\textit{Trip}_{peak}$ and $\textit{Trip}_\text{off-peak}$ denote the trips/person/h in peak and off-peak time, respectively. From~\cite{ABDELAAL2014677}, The average trip rate $\textit{Trip}_{mean}$ is 1.32 trips/person/day. Assume that there are only trips during peak times and off-peak times in a day. In~\cite{calabro2023adaptive}, peak hours are: 6h30 to 8h30 and 17h00 to 20h30 ($t_{peak}$ is 5.5 hours in total), and Off-peak hours are: 8h30 to 17h00 ($t_\text{off-peak}$ is 8.5 hours in total). And the ratio of average trip rate in one peak hour $\textit{Trip}_\text{peak}$ and one off-peak hour $\textit{Trip}_\text{off-peak}$ is 10:3. 
\begin{align}
&& \textit{Trip}_\text{peak} \cdot t_\text{peak} + \textit{Trip}_\text{off-peak} \cdot t_\text{off-peak} &= \textit{Trip}_\text{mean}
\label{eq:trips}
\end{align}

\begin{align}
&& \textit{Trip}_\text{peak} : \textit{Trip}_\text{off-peak} &= 10 : 3
\label{eq:trips-1}
\end{align}
Therefore, $\textit{Trip}_\text{peak} = 0.16 \,\,\,\, trips/person/hour$ is obtained by calculation. We replace this value of trip rate into Eqn.~\ref{eq:phi_i}-\ref{eq:flow-in}.

\subsection{Results}

\begin{figure}
    \centering
    \includegraphics[width=1.0\textwidth]{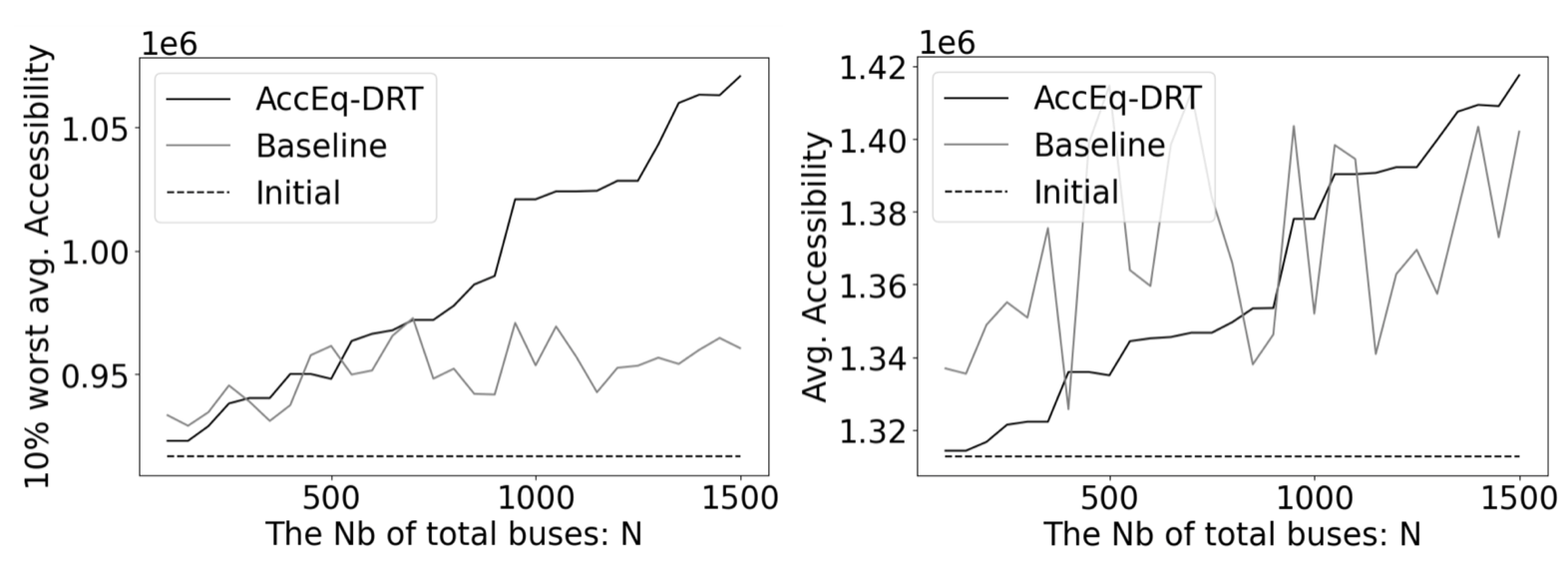}
    \caption{Improvement in accessibility obtained via DRT deployment } 
    \label{fig:Fig6}
\end{figure}

\begin{figure}
    \centering
    \includegraphics[width=0.5\textwidth]{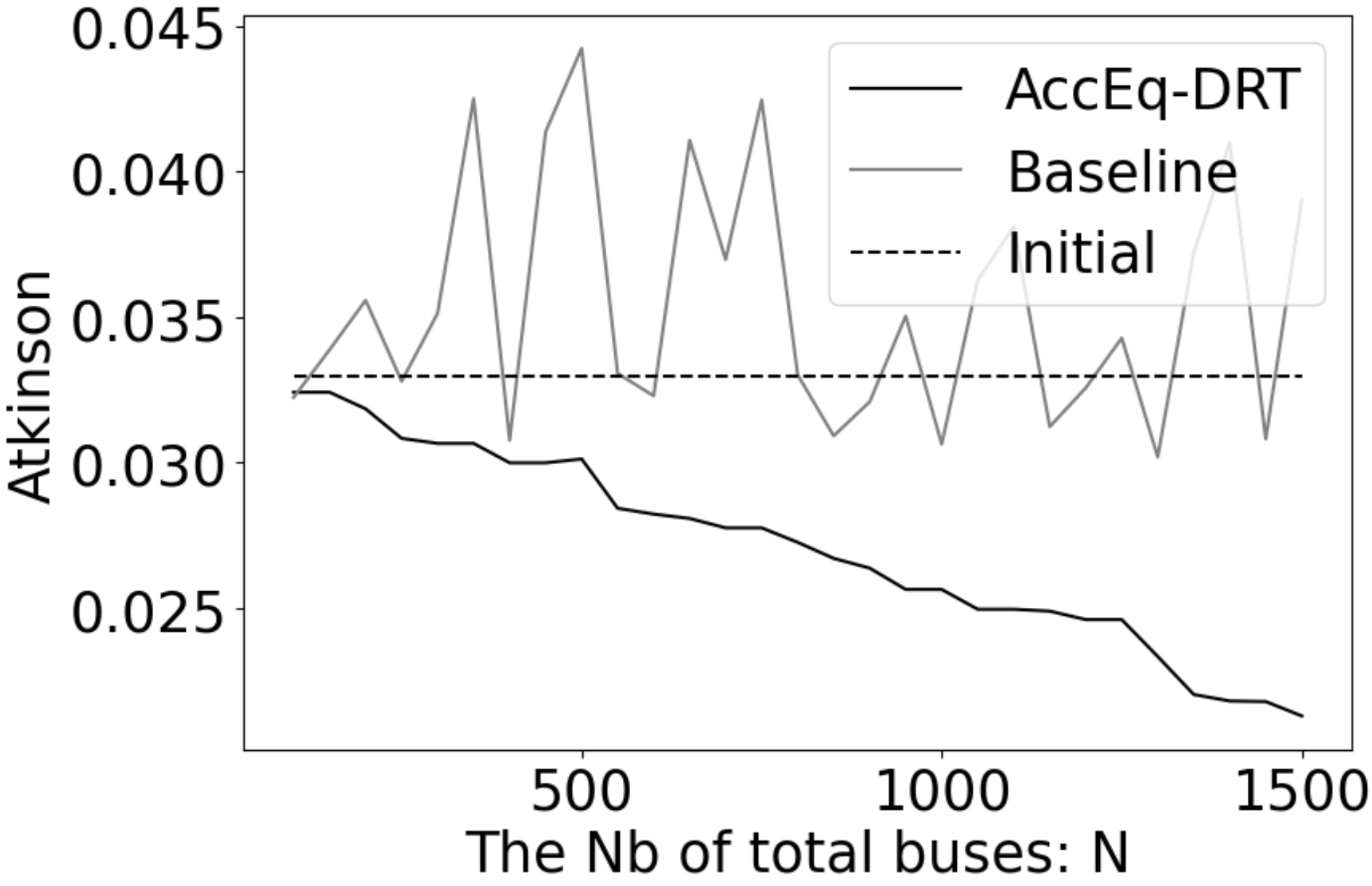}
    \caption{Inequality measured via the Atkinson index} 
    \label{fig:Fig7}
\end{figure}

We did not find in the literature any approach that combines a graph representing the real topology of a conventional PT network and an approximated CA model of DRT. We did not either find any approach to optimize DRT planning based on accessibility. A comparison with the state-of-the-art is thus necessarily non existent.

To show the relevance of the performance of AccEq-DRT, we compare it to a Baseline in which the number of DRT buses in each zone is chosen randomly. We ran the Baseline initially allowing to deploy DRT areas at potentially eveywhere. But then we realized that limiting the DRT zones in the Baseline, i.e., only serving with DRT the stops at the end of the lines, as in~\cite{steiner2020strategic}, was giving better results. Therefore, we keep this DRT zone limitation for the Baseline. 

The entire set of results concerning AccEq-DRT are obtained on an ordinary laptop in about 3h.

The values of fleet size range from 100 and 1500, which is a relatively conservative scenario, considering that in recent work (Fig.~5 of~\cite{Chouaki2021}) similar values have been considered for just a single suburb, while we are considering here an area equivalent to Montreal.

Fig.~\ref{fig:Fig6}-left shows that AccEq-DRT provides a strong improvement of the average accessibility of the individuals in the bottom decile, i.e., the 10\% individuals with the lowest accessibility. With 1500 DRT buses deployed in all Montreal, such individuals could reach more than 1.06 mln other people in a unit of time instead of 0.93 mln (16\% improvement). This suggests that AccEq-DRT reduces the gap of accessibility among individuals. This finding is confirmed by Fig.~\ref{fig:Fig7}, which shows that the Atkinson inequality index (Eqn.~\ref{eq:atkinson}) steadily decreases with the fleet size up to 34.4\% reduction with respect to conventional PT only. Fig.~\ref{fig:Fig7} also shows that deploying DRT without explicit consideration of inequality can exacerbate it (as the Baseline does). 
\begin{table}
	\caption{Performance in terms of different indexs}\label{tab:parameters_1}
	\begin{center}
		\begin{tabular}{l l l l l}
			Name & Initial ($a$) & AccEq-DRT ($b$)&  $\frac{a-b}{a}$ &  Reference\\\hline
			Atkinson & 0.032   & 0.021 & 34.4\%&  \cite{Costa2019inequality} \\
			Theil & 0.016 & 0.011 & 31.3\%&  \cite{Costa2019inequality} \\
            Pietra  & 0.074  & 0.058 & 21.6\%&  \cite{Costa2019inequality} \\
            Palma  & 0.400   & 0.371 & 7.3\%& \cite{sitthiyot2020simple}  \\\hline
		\end{tabular}
	\end{center}
\end{table}
The efficiency of AccEq-DRT is confirmed also when using inequality indicators different from Atkinson's (Table~\ref{tab:parameters_1}).

AccEq-DRT also improves overall accessibility and that this improvement steadily increases with the fleet size (Fig.~\ref{fig:Fig6}-right).

\begin{figure}
    \centering
    \includegraphics[width=0.7\textwidth]{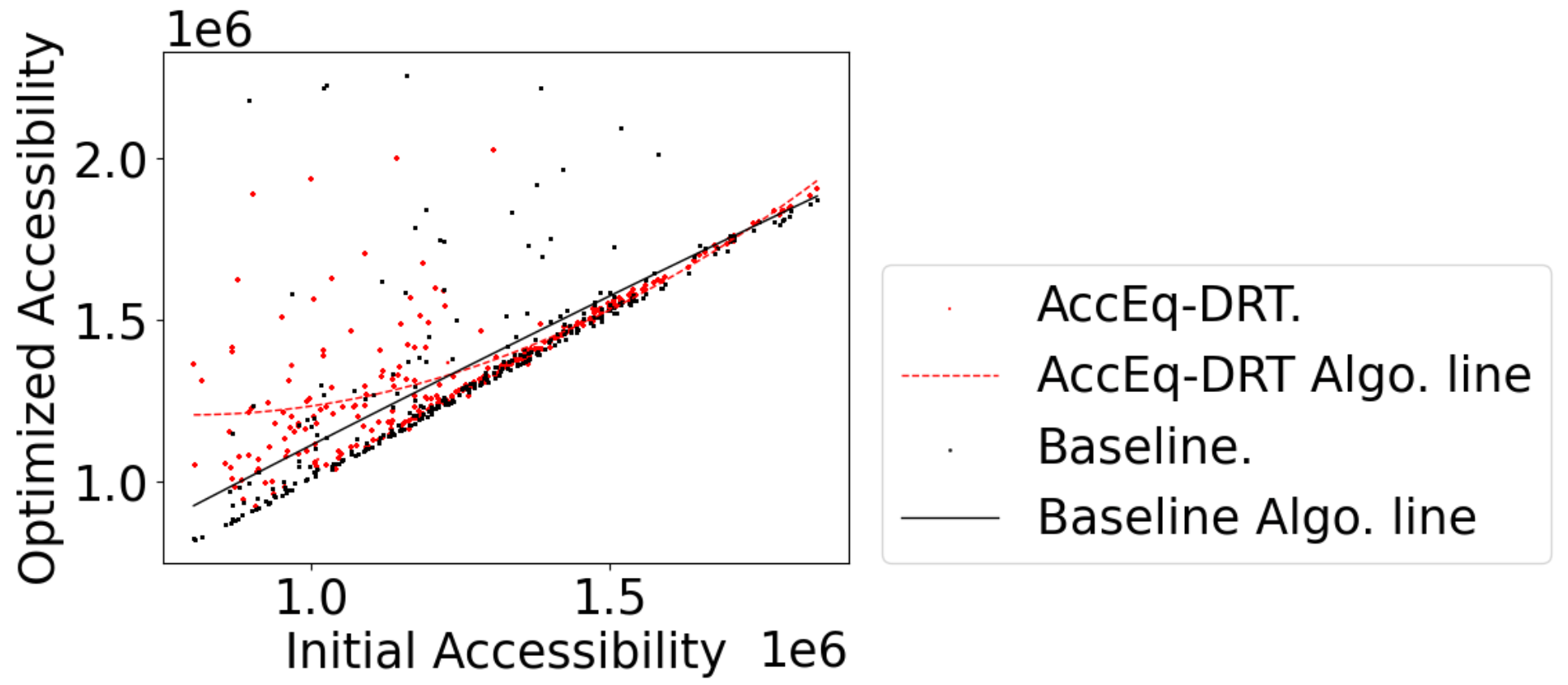}
    \caption{Distribution of accessibility of the two algorithms with the regression lines (AccEq-DRT vs. Baseline, with $N=1500$)} 
    \label{fig:Fig8}
\end{figure}

\begin{figure}
    \centering
    \includegraphics[width=0.65\textwidth]{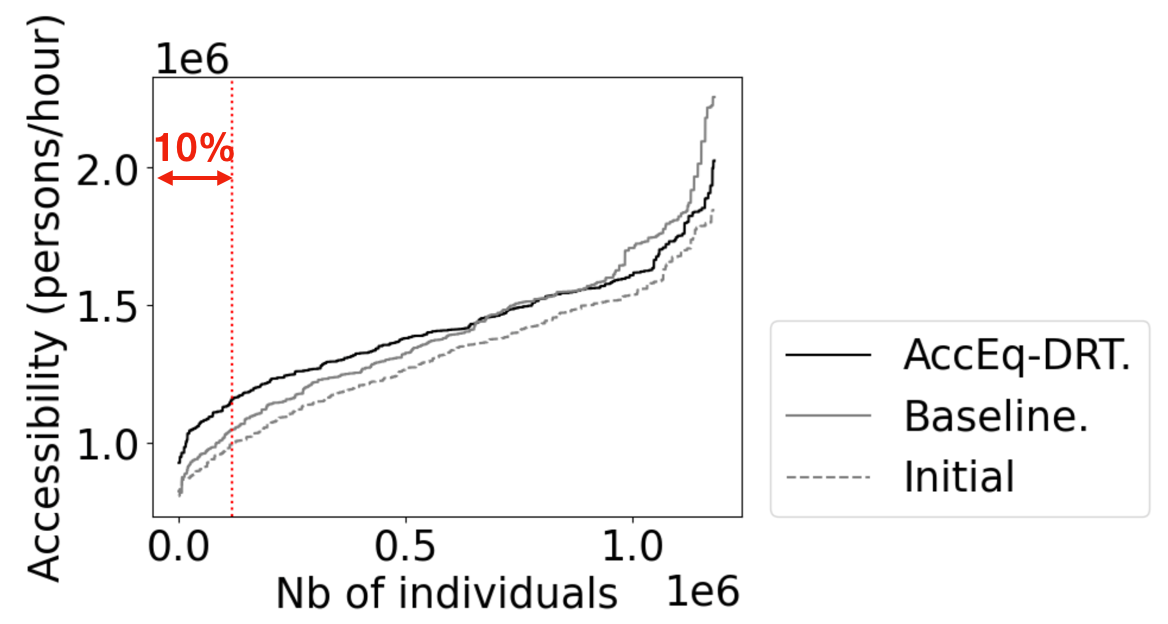}
    \caption{ Accessibility of each individual (with $N=1500$)} 
    \label{fig:Fig9}
\end{figure}
Fig.~\ref{fig:Fig8} and~\ref{fig:Fig9} explain the origin of the inequality reduction brought by AccEq-DRT. In these figures we compare the accessibility of each tile and each individual, before and after the introduction of DRT. As expected, introducing DRT always improves accessibility, no matter the deployment strategies. However, the Baseline distributes such improvement blindly across all tiles and individuals. AccEq-DRT, instead, improves the tiles and the individuals that were suffering from low accessibility with conventional PT only (those are the ones on the left of the x-axis of Fig.~\ref{fig:Fig8} and~\ref{fig:Fig9}). These results confirm that the score functions we designed in Eqn.~\ref{eq:index_c}-\ref{eq:index} are an effective guidance for the DRT bus allocation of Algorithm~\ref{fig:algo_2}.

\begin{figure}
    \centering
    \includegraphics[width=0.9\textwidth]{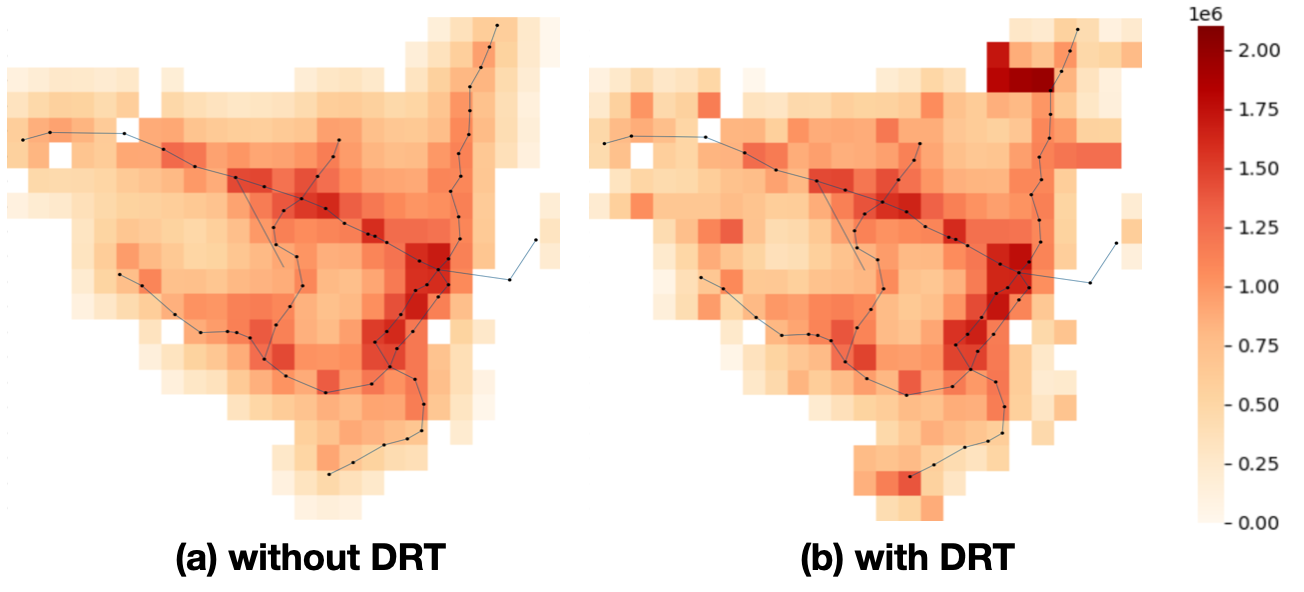}
    \caption{Heatmaps of accessibility in Montreal without DRT, with DRT ($N=1500$ buses). The white tiles are those removed from the study area (see Section ``Tessellation'').} 
    \label{fig:Fig5}
\end{figure}

\begin{figure}
    \centering
    \includegraphics[width=1.0\textwidth]{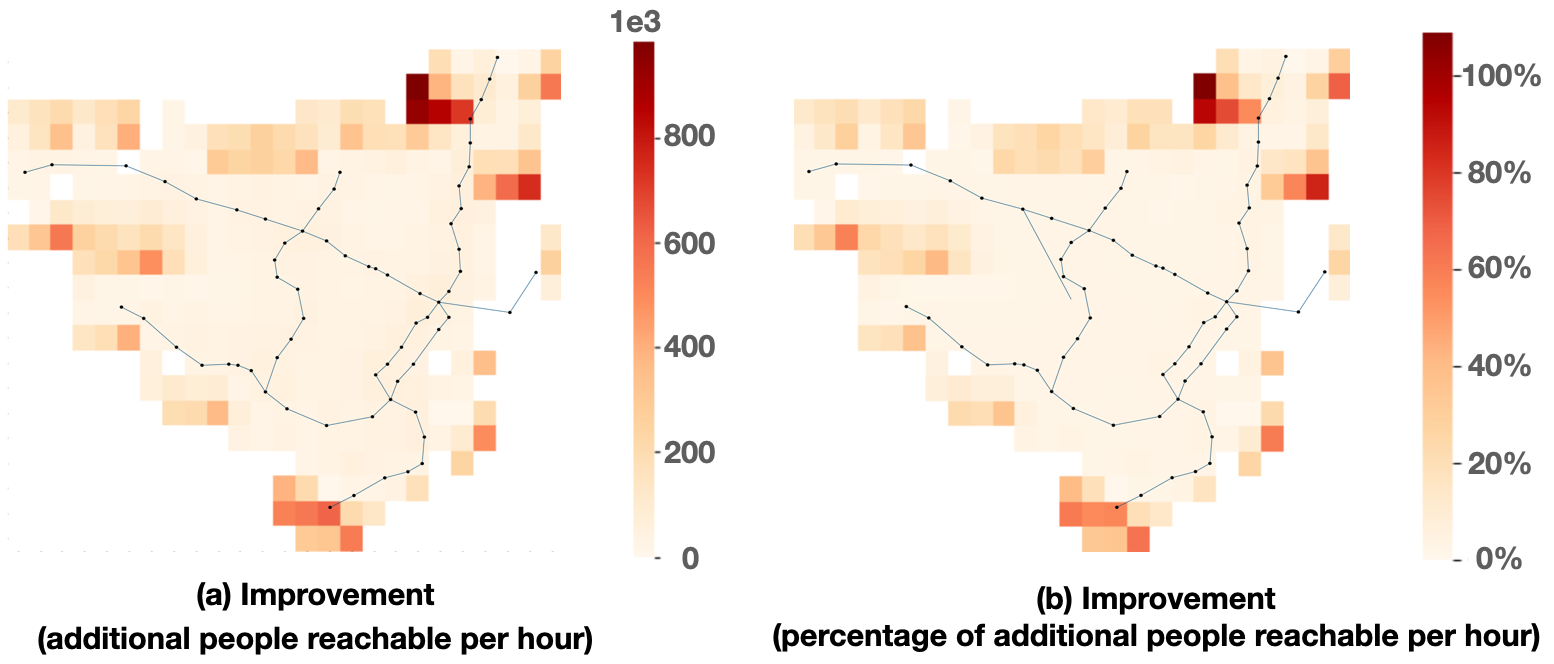}
    \caption{Heatmaps of the improvement achieved via AccEq-DRT ($N=1500$ buses)} 
    \label{fig:Fig13}
\end{figure}
Figure~\ref{fig:Fig5} shows the geographical distribution of accessibility in each tile before and after optimization by the AccEq-DRT model with 1500 DRT buses. In general, the accessibility of regions situated far from conventional PT service is relatively low before introducing DRT (Fig.~\ref{fig:Fig5}-(a)). In order to see the impact of AccEq-DRT more clearly, we calculate the improvement of accessibility for each tile.
For tile $c_i$, let us denote by~$acc^*(c_i)$ the accessibility of $c_i$ after optimization via the AccEq-DRT approach and by~$acc^0(c_i)$ the initial accessibility of $c_i$ (with conventional PT only).
Figure~\ref{fig:Fig13}-(a) depicts the improvement for tile $c_i$, i.e, $acc^*(c_i) - acc^0(c_i)$, which represents the additional people that are now reachable thanks to DRT, per unit of time. Figure~\ref{fig:Fig13}-(b) depictes the percentage improvement for tile $c_i$, i.e., $\frac{acc^*(c_i) - acc^0(c_i)}{acc^0(c_i)}$. We can see that AccEq-DRT effectively enhances accessibility of areas far from conventional PT (Fig.~\ref{fig:Fig13}). It is interesting to note that improvement is observed also where no DRT deployed. Although in such area is not possible to take a DRT bus, individuals residing there enjoy increase accessibility, as they can use DRT in the last mile to reach other people that were likely difficult to reach with conventional PT only.

\section{Conclusion}
We presented a methodology to plan DRT in order to reduce accessibility inequality among individuals in an urban region. We combine a graph model of conventional PT, capable of capturing the characteristics of currently existing PT, with a Continuous Approximation (CA) model of DRT. We combine the two modes in a single multilayer graph, on which we compute accessibility indicators. We propose a bilevel optimization method, where the upper level allocates buses in designated areas scattered around the city and the lower level solves transit assignment. Numerical results in an exemplary city shows strong improvement in the accessibility of those who were suffering low accessibility with conventional PT. Inequality is reduced.
 Our optimization is guided by a score that combines population and accessibility of the designated areas. While the structure of this score is fixed a-priori (based on our trial and error experiments), we will investigate the possibility to ``learn'' optimal indicators via Machine Learning (ML). By appropriately engineering features, ML would allow to consider complex local information (e.g. graph connectivity, topological or socio-economic data), which would be too difficult to capture with scoring functions constructed a-priori, and which could be customized for each city.

\section{CONTRIBUTIONS OF THE AUTHORS}

All authors contributed to the general concept and reviewed the paper. A.A. proposed the research goals, defined the methodological framework and provided funding. D.W. proposed AccEq-DRT optimization algorithm, implemented it and performed the analysis. M.E. provided advises and guidance for the model, the design of the algorithm and the analysis.

\section{Acknowledgements}
This work has been supported by the French ANR research project MuTAS (ANR-21-CE22-0025-01).

\bibliographystyle{unsrt}  
\bibliography{references}


\end{document}